\documentclass[reqno,11pt,a4]{amsart}
\usepackage{amsmath,amssymb,tabularx,setspace,color,enumerate}
\usepackage{multirow}

\usepackage{bigstrut}
\usepackage{multirow}
\usepackage{array}
\usepackage{tabularx,graphics}
\usepackage[pdftex]{graphicx}
\usepackage[margin=3.5cm]{geometry}

\newtheorem{definition}{Definition}
\newtheorem{remark}{Remark}

\makeatletter
\renewcommand\section{\@startsection {section}{1}{\z@}%
                                   {-3.5ex \@plus -1ex \@minus -.2ex}%
                                   {2.3ex \@plus.2ex}%
                                   {\normalfont\large\bfseries}}
\makeatother

\newtheorem{Theorem}{Theorem}
\begin{document}

\vspace{-0.2in}
\title[]{\vspace{-0.4in}Relationships between cumulative  entropy/extropy, Gini mean difference and probability weighted moments}

\author[]%
	{S\lowercase{udheesh} K. K\lowercase{attumannil}$\lowercase{^{a,\dag}}$,   S\lowercase{reedevi}, E. P.$\lowercase{^b}$  \lowercase{and} N. B\lowercase{alakrishnan}$\lowercase{{^c}}$  \\
$\lowercase{{^a}}$I\lowercase{ndian} S\lowercase{tatistical} I\lowercase{nstitute},
		C\lowercase{hennai}, I\lowercase{ndia},\\
			$\lowercase{{^b}}$M\lowercase{aharaja's} C\lowercase{ollege}, E\lowercase{ranakulam}, I\lowercase{ndia},\\
$\lowercase{{^c}}$M\lowercase{c}M\lowercase{aster} U\lowercase{niversity}, H\lowercase{amilton}, C\lowercase{anada}.}
	\thanks{{$^{\dag}$Corresponding author E-mail: \tt skkattu@isichennai.res.in}.}
	
\vspace{-0.8in}
\doublespace
\begin{abstract}
 In this work, we establish a connection between the cumulative residual entropy and the Gini mean difference (GMD).  Some relationships between the extropy and the GMD, and the truncated GMD and dynamic versions of the cumulative extropy are also established.
 We then show that several entropy and extropy measures discussed here can be brought into the framework of probability weighted moments, which would enable us to find estimators of these measures.
\end{abstract}
\maketitle
\vspace{-0.1in}
\section{Introduction}\vspace{-0.1in}The  Gini mean difference (GMD) is a prominent measure that is used extensively in economics.
 Let $X$ be a non-negative random variable with absolutely  continuous distribution function $F$  and   finite mean $\mu$, and   $X_1$ and $X_2$ be two  independent random variables from $F$.  Then, the GMD is defined as
\begin{equation*}
  GMD=E|X_1-X_2|.
\end{equation*}
For more details on  GMD and other measures derived from it, interested readers may refer to Yitzhaki and Schechtman (2013).

The truncated version of  the GMD has also been discussed for studying inequality prevailing in the poor and  affluent groups,  while the left truncated version has found application in reliability analysis.  For more details, see Nair et al. (2012) and  Behdani et al. (2020).   Recently, Nair and Vineshkumar (2022) discussed some relationships between the cumulative entropy  and income gap ratio, Lorenz curve, Gini index, Bonferroni curve and  Zenga Curve.

 In many practical situations, measuring the uncertainty associated with a random variable is quite  important, and  many measures have been introduced for this purpose.  The seminal work on information theory started with the concept of Shannon entropy or differential entropy introduced by Shannon (1948). Since then, different measures of entropy have been discussed, each one being suitable for some specific situations. A few of the widely used measures of entropy are cumulative residual entropy (Rao et al., 2004), cumulative entropy (Di Crescenzo and Longobardi, 2009), the corresponding weighed measures by Mirali et al. (2016) and Mirali and Baratpour (2017), and a general measure of cumulative residual entropy (Sudheesh et al., 2022).

One important generalization of Shannon entropy is due to Tsallis (1988), known as generalized Tsallis entropy of order $\alpha$.
Many extensions and modifications are available for it as well. For example,  Rajesh and Sunoj (2019) proposed cumulative residual Tsallis entropy of order $\alpha$, while  Chakraborty and Pradhan (2021a) defined weighted cumulative residual Tsallis  entropy (WCRTE) of order $\alpha$ and its dynamic version.
 Calì et al. (2017) introduced cumulative Tsallis past entropy of order $\alpha$. Chakraborty and Pradhan (2021b) introduced  weighted cumulative Tsallis entropy (WCTE) of order $\alpha$,  and also studied its dynamic version.

An alternative measure of uncertainty, called extropy, was introduced by Lad et al. (2015) as a complementary dual of entropy.
Jahanshahi et al. (2020) and Tahmasebi and Toomaj (2022) studied cumulative residual extropy and negative cumulative extropy, while   Balakrishnan et al. (2022) and  Chakraborty  and  Pradhan  (2022) discussed different weighted versions of extropy.  Sudheesh and Sreedevi (2022) established some   relationships between different extropy measures  and reliability concepts, and Sudheesh et al. (2022) established some relationships between entropy and extropy measures.

The probability weighted moments (PWM) generalize the concept of moments of a probability distribution, and they have been used effectively for estimating the parameters of different distributions.  The PWM was introduced by Greenwood et al. (1979).   Here, we show that may  extropy measures mentioned  above can be expressed in terms of PWM.

The rest of this article is organised as follows.   In Section 2,  using a  generalized cumulative residual entropy, we derive connections between some entropy measures and GMD. In Section 3, we establish some relationships between GMD and extropy measures. In Section 4, we show that several extropy measures can be brought into the framework of PWM. Finally,  we make some concluding remarks in Section 5.

\section{Connection between entropy and GMD}
 Let $X$ be a non-negative  random variable with absolute continuous distribution function  $F$ and let $\bar{F}(x)=P(X>x)$ denote the survival function. Let us further assume that the mean $\mu=E(X)<\infty$. 
\subsection{Cumulative residual entropy and GMD}
We now consider  the generalized cumulative residual entropy measure introduced recently by Sudheesh et al. (2022).
\begin{definition}
  Let $X$ be a non-negative  random variable with absolute continuous distribution function  $F$. Further, let $ \phi(.)$ be  a function  of $X$ and $w(.)$ be a weight function.  Then, the generalized cumulative residual entropy is defined  as
\begin{equation}\label{gwentropy}
  \mathcal{GE}(X)=\int w(u) E\Big[\phi(X)-\phi(u)| X> u\Big]dF(u).
\end{equation}
\end{definition}

  Sudheesh et al. (2022) then showed that several measures of entropy available in the literature can be deduced from (\ref{gwentropy}) with different choices of  $w(\cdot)$ and $\phi(\cdot)$.
  We now show that GMD is a special case of  $\mathcal{GE}(X)$.\vspace{-0.2in}
  \begin{Theorem}
    For the choice of $w(u)=\bar{F}(u)$ and $\phi(x)=2x$, $\mathcal{GE}(X)=GMD$.
  \end{Theorem}\vspace{-0.2in}
  \noindent{\bf Proof:}
  Let $X_1$ and $X_2$ be independent and identical  random variables having the distribution function $F$. With $  GMD=E|X_1-X_2|$, and because $|X_1-X_2|=\max(X_1,X_2)-\min(X_1,X_2)$, the GMD can alternatively be expressed as\vspace{-0.1in}
\begin{equation}\label{gmd}
  GMD=E\left(\max(X_1,X_2)-\min(X_1,X_2)\right).
\end{equation}\vspace{-0.2in}
Now, for the choice of $w(u)=\bar{F}(u)$ and $\phi(x)=2x$, the  expression in (\ref{gwentropy}) becomes
\begin{eqnarray*}
  \mathcal{GE}(X)&=&2\int_{0}^{\infty}\bar{F}(u)E(X-u|X>u)dF(u)\\
   &=& 2\int_{0}^{\infty}\bar{F}(u)\left(\frac{1}{\bar{F}(u)}\int_{u}^{\infty}(y-u)dF(y)\right)dF(u)\\
   &=&2\int_{0}^{\infty}\int_{u}^{\infty}ydF(y)dF(u)-2\int_{0}^{\infty}u\int_{u}^{\infty}dF(y)dF(u)
   \\
   &=&\int_{0}^{\infty}y\cdot2F(y)dF(y)-\int_{0}^{\infty}y\cdot2\bar{F}(y)dF(y)\\
   &=&E(\max(X_1,X_2))-E(\min(X_1,X_2)),
\end{eqnarray*}which is precisely the GMD.

Instead,  if we choose  $w(u)=\bar{F}(u)$ and $\phi(x)=2x^2$,  (\ref{gwentropy}) becomes
\begin{eqnarray*}
  \mathcal{GE}(X) &=&2\int_{0}^{\infty}\int_{u}^{\infty}y^2dF(y)dF(u)-2\int_{0}^{\infty}u\int_{u}^{\infty}dF(y)dF(u)
   \\
   &=&\int_{0}^{\infty}y^2\cdot 2 F(y)dF(y)-\int_{0}^{\infty}y^2\cdot 2\bar{F}(y)dF(y)\\
   &=&E(\max(X_1,X_2)^2)-E(\min(X_1,X_2)^2).
\end{eqnarray*}In general, for the choices of $w(u)=\bar{F}(u)$ and  $\phi(x)=2x^v$ for some positive integer $v$,  we similarly have\vspace{-0.2in}
\begin{eqnarray*}
  \mathcal{GE}(X)=E\left(\max(X_1,X_2)^v-\min(X_1,X_2)^v\right).
\end{eqnarray*}
Also, for  $w(x)=\bar{F}^{k-1}(u)$, $k=2,3,\ldots$ and $\phi(x)=kx^v$,  we get
\begin{equation*}
   \mathcal{GE}(X)=E\left(\max(X_1,\ldots,X_{k})^v-\min(X_1,\ldots,X_{k})^v\right).
\end{equation*}
In the particular case of  $v=1$, we similarly  have
\begin{eqnarray*}
   \mathcal{GE}(X)&=&E(\max(X_1,\ldots,X_{k})-\min(X_1,\ldots,X_{k})) \\
   &=& E(\max(X_1,\ldots,X_{k}))-E(X_1)+E(X_1)-(E(\min(X_1,\ldots,X_{k})) \\
   &=&EG_k(-X)-EG_k(X),
\end{eqnarray*}where $EG_k(X)$ and $EG_k(-X)$ are the risk-premium and the gain-premium of $X$ of
order $k$, respectively. These two quantities have found key applications in bid and ask prices in finance; see Agouram and Lakhnati (2016).
The above expression can alternatively be expressed as
\begin{eqnarray*}
   \mathcal{GE}(X)&=& E(\max(X_1,\ldots,X_{k}))-E(X_1)+E(X_1)-(E(\min(X_1,\ldots,X_{k})) \\
   &=&kCov(X,F^{k-1}(X))-kCov(X,\bar F^{k-1}(X)).
\end{eqnarray*}\vspace{-0.1in}In particular, when $k=2$, we deduce
\vspace{-0.1in}
   $$\mathcal{GE}(X)=4Cov(X,F(X)),$$which is yet another representation of GMD.

\subsection{Generalized  cumulative entropy and GMD }
We now consider  the generalized cumulative  entropy and discuss its relationship with GMD.
\begin{definition}[Sudheesh et al., 2022]
  Let $X$ be a non-negative  random variable with absolute continuous distribution function  $F$. Further,  let $ \phi(\cdot)$ be  a function  of $X$ and $w(\cdot)$ be a weight function.  Then, the generalized cumulative  entropy is defined  as
\begin{equation}\label{gentropy}
  \mathcal{GCE}(X)=\int_{0}^{\infty}w(u)E\Big[\phi(u)-\phi(X)|X\le u\Big]dF(u).
\end{equation}
\end{definition}
In particular, for the choice of  $w(u)=1$ and $\phi(x)=x$,  (\ref{gentropy}) reduces to $\mathcal{CE} (X)$ (Di Crescenzo and Longobardi, 2009).
   \begin{Theorem}
      For the choice of $w(u)=F(u)$ and $\phi(x)=x$, $\mathcal{GCE}(X)$ reduces to the GMD.
   \end{Theorem}\vspace{-0.2in}
   \noindent{\bf Proof:} Consider\vspace{-0.1in}
     \begin{eqnarray*}\vspace{-0.2in}
  \mathcal{GE}(X)&=&2\int_{0}^{\infty}{F}(u)E(u-X|X<u)dF(u)\\
   &=& 2\int_{0}^{\infty}{F}(u)\left(\frac{1}{{F}(u)}\int_{0}^{u}(u-y)dF(y)\right)dF(u)\\
   &=&2\int_{0}^{\infty}\int_{0}^{u}udF(y)dF(u)-2\int_{0}^{\infty}\int_{0}^{u}ydF(y)dF(u)
   \\
   &=&E\left(\max(X_1,X_2)\right)-\int_{0}^{\infty}y\cdot 2\bar{F}(y)dF(y)\\
   &=&E\left(\max(X_1,X_2)-\min(X_1,X_2)\right).
\end{eqnarray*}\vspace{-0.15in} Further, for the choice of $w(u)=F^k(u)$ and $\phi(x)=(k+1)x$, we have
\begin{equation*}
  \mathcal{GCE}(X)=E\left(\max(X_1,\ldots,X_{k+1})-\min(X_1,\ldots,X_{k+1})\right).
\end{equation*}

\vspace{-0.1in} We now show that GMD is a special case of the cumulative residual Tsallis entropy of order $\alpha$ defined by Rajesh and Sunoj (2019) as
\begin{equation}\label{tsent}
  CRT_{\alpha}(X)=\frac{1}{\alpha-1}\int_{0}^{\infty} (\bar{F}(x)-\bar F^{\alpha}(x))dx, \,\alpha>0,\,\alpha\ne 1.
\end{equation}For $\alpha=2$, the above expression becomes
\begin{equation}\label{tsent}\vspace{-0.1in}
  CRT_{2}(X)=\int_{0}^{\infty} (\bar{F}(x)-\bar F^{2}(x))dx.
\end{equation}For a non-negative random variable $X$, we have $\mu=E(X)=\int_{0}^{\infty}\bar F(x)dx$.  Noting that $\bar F^{2}(x)$ is the survival function of $\min(X_1,X_2)$, we simply obtain
\begin{equation}\label{tsent}\vspace{-0.1in}
  CRT_{2}(X)=E(X_1)-E(\min(X_1,X_2)).
\end{equation} Now, because $|X_1-X_2|=X_1+X_2-2\min(X_1,X_2)$, we obtain $GMD=2 CRT_{2}(X)$.\vspace{-0.2in}
\section{Connection between extropy and GMD}\vspace{-0.1in} Tahmasebi and Toomaj (2022) discussed a relationship between GMD and cumulative extropy. Here using the generalized  entropy measures introduced by Sudheesh et al. (2022), we establish some relationships between different dynamic extropy measures and the right and left truncated GMD.

First, we express the GMD  in terms of cumulative residual extropy and residual extropy.  For a non-negative random variable  $X$, Jahanshahi et al. (2020) defined the cumulative residual extropy as
\begin{equation}\label{creext}
   \mathcal{CRJ}(X)=-\frac{1}{2}\int_{0}^{\infty}\bar{F}^2(x)dx,
\end{equation}while the cumulative  extropy is defined as (Tahmasebi and Toomaj,  2020)
\begin{equation}\label{cext}
   \mathcal{CJ}(X)=-\frac{1}{2}\int_{0}^{\infty}(1-{F}^2(x))dx.
\end{equation}Using the survival functions of $\max(X_1,X_2)$ and $\min(X_1,X_2)$, we easily find that
\begin{eqnarray*}
 2 \mathcal{CRJ}(X)-2\mathcal{CJ}(X)&=&\int_{0}^{\infty}(1-{F}^2(x))dx-\int_{0}^{\infty}\bar{F}^2(x)dx\\
  &=&E\left(\max(X_1,X_2)-\min(X_1,X_2)\right)\\
  &=&GMD.
\end{eqnarray*}Sudheesh and Sreedevi (2022) discussed the non-parametric estimation of $\mathcal{CRJ}(X)$ and $\mathcal{CJ}(X)$ for right censored data,  and so using the above relationship, we can readily obtain an estimator of the GMD based on right censored data.

Sathar and Nair (2021) defined  dynamic survival extropy as
\begin{equation}\label{dsurvent}
  J_t(X)=-\frac{1}{2\bar{F}^2(t)}\int_{t}^{\infty}\bar{F}^2(x)dx,
\end{equation}
 while Kundu (2021) introduced dynamic cumulative extropy as
\begin{equation}\label{dcext}
  H_t(X)=-\frac{1}{2{F}^2(t)}\int_{0}^{t}{F}^2(x)dx.
\end{equation} Sudheesh and Sreedevi (2022) proposed simple alternative expressions for different extropy measures, and then  used them to establish some relationships between different dynamic and weighted extropy measures and reliability concepts.

 The left truncated GMD is given by
\begin{equation*}
  GMD_L=\frac{1}{\bar{F}^2(t)}\int_{t}^{\infty} (\bar{F}(t)-2\bar{F}(x))xdF(x)
\end{equation*}while the right truncated GMD is defined as
\begin{equation*}
  GMD_R=\frac{1}{{F}^2(t)}\int_{0}^{t} (2{F}(x)-{F}(t))xdF(x).
\end{equation*}
\begin{Theorem}Let $m(x)=E(X-x|X>x)$ and $r(x)=E(x-X|X\le x)$. Then, the  relationships between the GMD and  the dynamic survival extropy and/or the dynamic cumulative extropy are \\
(i) \begin{equation}\label{gmdlent}
 GMD_L=m(x)+2J_t(X),
\end{equation}
(ii) \begin{equation}\label{gmdcent}
  GMD_R=-2H_t(x)-r(x).
\end{equation}
\end{Theorem}

\noindent{\bf Proof:} Consider the left truncated GMD given by
\begin{eqnarray}\label{lgmdn}
  GMD_L&=&\frac{1}{\bar{F}^2(t)}\int_{t}^{\infty} (\bar{F}(t)-2\bar{F}(x))xdF(x)\nonumber\\
  &=&\frac{1}{\bar{F}^2(t)}\int_{t}^{\infty} \left(\bar{F}(t)-t+t-2\bar{F}(x)\right)xdF(x)\nonumber\\
  &=&\frac{1}{\bar{F}(t)}\int_{t}^{\infty} (x-t)dF(x)-\frac{1}{\bar{F}^2(t)}\int_{t}^{\infty} (x-t)2\bar{F}(x)dF(x)\nonumber\\
  &=&m(x)-E(\min(X_1,X_2)-t|\min(X_1,X_2)>t).
\end{eqnarray}Sudheesh and Sreedevi (2022) expressed $J_t(X)$ in (\ref{dsurvent}) as
\begin{equation}\label{dsextropy}
  J_t(X)=-\frac{1}{2}E\left(\min(X_1,X_2)-t|\min(X_1,X_2)>t\right).
\end{equation}
Substituting  (\ref{dsextropy}) in (\ref{lgmdn}), we obtain the relationship in (\ref{gmdlent}).

Next, consider the right truncated GMD given by
\begin{eqnarray}\label{rgmdn}
  GMD_R&=&\frac{1}{{F}^2(t)}\int_{0}^{t} (2{F}(x)-{F}(t))xdF(x).
\end{eqnarray}After some algebraic manipulation, we obtain from (\ref{rgmdn}) that
\begin{eqnarray}\label{rgmdnn}
  GMD_R=E\left(t-\min(X_1,X_2)|\min(X_1,X_2)\le t\right)-r(x).
\end{eqnarray}Sudheesh and Sreedevi (2022) expressed $H_t(X)$ as
\begin{equation}\label{dcextropy}
H_t(X)=-\frac{1}{2}E\left(t-\max(X_1,X_2)|\max(X_1,X_2)\le t\right).
\end{equation}
Substituting  (\ref{dcextropy}) in (\ref{rgmdnn}) we obtain the relationship in (\ref{gmdcent}).

\begin{Theorem}
  Let $X_1$ and $X_2$ be two independent and identical  random variables having distribution function $F$, and $Z=\min(X_1,X_2)$ be the lifetime of a  series system with the two components. Then, the generalized residual entropy associated with $Z$ is the weighted average of the difference between the left truncated  GMD and the mean residual life.
\end{Theorem}
\noindent{\bf Proof:}    Using (\ref{gentropy}), we define the generalized residual entropy associated with $Z$ as
\begin{equation}\label{gwentropyz}
  \mathcal{GE}(Z)=\int 2w(u) E\Big[\phi(Z)-\phi(u)| Z> u\Big]\bar{F}(u)dF(u),
\end{equation}where $ \phi(.)$ is  a function  of $Z$ and $w(.)$ is a weight function.
Now, for the choice of $w(u)=-\frac{1}{2}$ and $\phi(z)=z$,  from (\ref{gwentropyz}), we obtain
\begin{equation*}\label{gwentropyzn}
  \mathcal{GE}(Z)=\int (GMD_L-m(u))\bar{F}(u)dF(u).
\end{equation*}
\begin{Theorem}Let $X_1$ and $X_2$ be two independent and identical  random variables having distribution function $F$, and $Z=\max(X_1,X_2)$ be the lifetime of a  parallel  system with the two  components. Then, the generalized cumulative  entropy associated with $Z$ is the weighted average  of the  sum of the right truncated  GMD and the mean past life.
\end{Theorem}
\noindent{\bf Proof:}
The generalized cumulative  entropy associated with $Z$ is defined  as
\begin{equation}\label{gentropyz}
  \mathcal{GCE}(Z)=\int_{0}^{\infty}2w(u)E\Big[\phi(u)-\phi(X)|X\le u\Big]F(u)dF(u).
\end{equation}
Again,  for the choice of $w(u)=\frac{1}{2}$ and $\phi(z)=z$,  from (\ref{gentropyz}), we obtain
\begin{equation*}\label{gentropyzn}
  \mathcal{GCE}(Z)=-\int_{0}^{\infty}(GMD_R+r(x))F(u)dF(u).
\end{equation*}So, the generalized cumulative entropy associated with $Z$ is  the weighted average of sum of the right truncated  GMD and the mean past life.

\section{Connections with probability weighted moments}
We now show that several measures discussed in the preceding sections can be represented in terms of PWM.

The PWM of a random variable $X$  with distribution function  $F$ is defined as (Greenwood et al., 1979)
\begin{equation}\label{pwm}
\mathcal{M}_{p,r,s}=\mathbb{E}\left \lbrace X^pF^r(X)(1-F(X))^s\right\rbrace,
\end{equation}
where $p$, $r$ and $s$ are any real numbers, for which the involved expectation exists.

First, we consider the GMD and related measures.   From the proof of Theorem 1, we can express GMD as
\begin{equation*}
  GMD=\int_{0}^{\infty}2yF(y)dF(y)-\int_{0}^{\infty}2y\bar{F}(y)dF(y),
\end{equation*} and so we have
\begin{equation*}
  GMD=2E(XF(X))-2E(X\bar{F}(X)),
\end{equation*}
which, when compared with (\ref{pwm}), yields
\begin{equation*}
  GMD=2\mathcal{M}_{1,1,0}-2\mathcal{M}_{1,0,1}.
\end{equation*}Several income inequality measures are derived from the GMD by choosing different
weights in the expectation and one among them is the S-Gini family of indices (Yitzhaki and Schechtman, 2013).  The absolute S-Gini index is defined as
\begin{equation*}
  S_v=-Cov(X,\bar{F}^{v-1}(X)),\quad v>0 \, \text{ and }\, v\ne 1,
\end{equation*} which  can be rewritten  as
\begin{equation*}
  S_v=-E(X\bar{F}^{v-1}(X))+\frac{1}{v}E(X).
\end{equation*}  Hence, the absolute S-Gini index can be represented as
\begin{equation*}
 S_v=\frac{1}{v}\mathcal{M}_{1,0,0}-\mathcal{M}_{1,0,v-1}.
\end{equation*}

Next, we consider different extropy measures. The cumulative extropy and the  weighted cumulative extropy can be expressed as (Sudheesh and Sreedevi, 2022)
\begin{equation*}\label{adef}
  \mathcal{CRJ}(X)=-\frac{1}{2}E(\max(X_1,X_2)),
\end{equation*}
\begin{equation*}\label{wadef}
  \mathcal{CRJW}(X)=-\frac{1}{4}E(\max(X_1,X_2)^2).
\end{equation*}
Observing that  $2F(x)f(x)$ is the density function of $\max(X_1,X_2)$,  the extropy measures given above can be rewritten as
\begin{equation*}\label{cdef1}
  \mathcal{CRJ}(X)=-E(XF(X)),
\end{equation*}
\begin{equation*}\label{wcdef1}
  \mathcal{CRJW}(X)=-\frac{1}{2}E(X^2F(X)),
\end{equation*}
so that
\begin{equation*}
  \mathcal{CRJ}(X)=-\mathcal{M}_{1,1,0}\quad \text{ and }\quad  \mathcal{CRJW}(X)=-\frac{1}{2}\mathcal{M}_{2,1,0}.
\end{equation*}

The cumulative residual entropy and the weighted  cumulative residual entropy can be expressed as (Sudheesh and Sreedevi, 2022)
\begin{equation*}
  \mathcal{CE}(X)=-\frac{1}{2}E(\min(X_1,X_2)),
\end{equation*}
\begin{equation*}
  \mathcal{WCE}(X)=-\frac{1}{4}E(\min(X_1,X_2)^2).
\end{equation*}
Here again, we observe that  $2\bar{F}(x)f(x)$ is the density function of $\min(X_1,X_2)$, and so
\begin{equation*}\label{edef1}
  \mathcal{CE}(X)=-E(X\bar{F}(X))
\end{equation*}
\begin{equation*}\label{wedef1}
  \mathcal{WCE}(X)=-\frac{1}{2}E(X^2\bar{F}(X))),
\end{equation*}
from which we readily obtain
\begin{equation*}
  \mathcal{CE}(X)=-\mathcal{M}_{1,0,1}\quad \text{ and }\quad  \mathcal{WCE}(X)=-\frac{1}{2}\mathcal{M}_{2,0,1}.
\end{equation*}

Next, we show that the cumulative (residual) Tsallis entropy of order $\alpha$ can be written in terms of PWM. Let us assume that $\alpha$ is an integer.    Recall that
\begin{equation}\label{tsentn}
  CRT_{\alpha}(X)=\frac{1}{\alpha-1}\int_{0}^{\infty} (\bar{F}(x)-\bar F^{\alpha}(x))dx.
\end{equation}
We can show that
\begin{eqnarray*}
 \int_{0}^{\infty}x \alpha \bar F^{\alpha-1}(x)dF(x) =\int_{0}^{\infty}\bar F^{\alpha}(x)dx,
\end{eqnarray*} and so we can write  from (\ref{tsentn}) that
\begin{equation*}\label{tsentnp}
  CRT_{\alpha}(X)=\frac{1}{\alpha-1}E(X)-\frac{\alpha}{\alpha-1}E(X\bar F^{\alpha}(X)),
\end{equation*}which yields
\begin{equation*}
   CRT_{\alpha}(X)=\frac{1}{\alpha-1} \mathcal{M}_{1,0,0}-\frac{\alpha}{\alpha-1}\mathcal{M}_{1,0,(\alpha-1)}.
\end{equation*}

For a non-negative continuous random variable $X$,  the weighted cumulative residual Tsallis  entropy (WCRTE) of order $\alpha$ is defined as (Chakraborty and Pradhan, 2021a)
\begin{equation*}\label{wctre}
  WCRT_{\alpha}(X)=\frac{1}{\alpha-1}\int_{0}^{\infty}x(\bar{F}(x)-\bar{F}^{\alpha}(x))dx.
\end{equation*}After some simple algebra, we can show that
\begin{eqnarray*}
 \int_{0}^{\infty}x^2 \alpha \bar F^{\alpha-1}(x)dF(x) =\int_{0}^{\infty}2x\bar F^{\alpha}(x)dx,
\end{eqnarray*}so that
\begin{equation*}\label{wctren}
  WCRT_{\alpha}(X)=\frac{1}{2(\alpha-1)}E(X^2)-\frac{\alpha}{2(\alpha-1)}E(X^2\bar{F}^{\alpha}(X))
\end{equation*} from which we readily obtain
\begin{equation*}
   WCRT_{\alpha}(X)=\frac{1}{2(\alpha-1)} \mathcal{M}_{2,0,0}-\frac{\alpha}{2(\alpha-1)}\mathcal{M}_{2,0,(\alpha-1)}.
\end{equation*}
Calì et al. (2017) introduced the cumulative Tsallis past entropy as
\begin{equation*}
   CT_{\alpha}(X)=\frac{1}{\alpha-1}\int_{0}^{\infty} ({F}(x)- F^{\alpha}(x))dx,\,\alpha>0, \alpha\ne1,
\end{equation*}
while Chakraborty and Pradhan (2021b) defined the weighted  cumulative Tsallis entropy (WCTE) of order $\alpha$ as
\begin{equation}\label{wcte}
  WCT_{\alpha}(X)=\frac{1}{\alpha-1}\int_{0}^{\infty}x(F(x)-F^{\alpha}(x))dx,\,\alpha>0, \alpha\ne1.
\end{equation}
As in the case of CRTE and WCRTE, we obtain
\begin{equation*}
   CT_{\alpha}(X)=\frac{1}{(\alpha-1)} \mathcal{M}_{1,0,0}-\frac{\alpha}{(\alpha-1)}\mathcal{M}_{1,(\alpha-1),0},
\end{equation*}
\begin{equation*}
   WCRT_{\alpha}(X)=\frac{1}{2(\alpha-1)} \mathcal{M}_{2,0,0}-\frac{\alpha}{2(\alpha-1)}\mathcal{M}_{2,(\alpha-1),0}.
\end{equation*}
Sharma and Taneja (1975) and Mittal (1975) independently introduced the STM entropy  as
\begin{equation}\label{stm}
  S_{\alpha,\beta}=\frac{1}{\beta-\alpha}\int_{0}^{\infty}(f^{\alpha}(x)-f^{\beta}(x))dx.
\end{equation}
Some entropy measures discussed in the literature can be derived from (\ref{stm}) by  taking different values  of  $\alpha$ and $\beta$.

Sudheesh et al. (2022) introduced a  cumulative residual STM  entropy as
\begin{equation*}\label{stmr}
  SR_{\alpha,\beta}=\frac{1}{\beta-\alpha}\int_{0}^{\infty}(\bar F^{\alpha}(x)-\bar F^{\beta}(x))dx,
\end{equation*}
and the cumulative   past  STM  entropy  as
\begin{equation*}\label{stmp}
  SP_{\alpha,\beta}=\frac{1}{\beta-\alpha}\int_{0}^{\infty}(F^{\alpha}(x)-F^{\beta}(x))dx.
\end{equation*}
They  also introduced weighted versions of these measures as
\begin{equation*}\label{stmr}
  SRW_{\alpha,\beta}=\frac{1}{\beta-\alpha}\int_{0}^{\infty}x(\bar F^{\alpha}(x)-\bar F^{\beta}(x))dx
\end{equation*}
 and
\begin{equation*}\label{stmpw}
  SPW_{\alpha,\beta}=\frac{1}{\beta-\alpha}\int_{0}^{\infty}x(F^{\alpha}(x)-F^{\beta}(x))dx.
\end{equation*}
We now express these measures in terms of PWM as follows, presented without proofs, for the sake of conciseness:\vspace{-0.2in}
\begin{equation*}\label{stmr}
  SR_{\alpha,\beta}=\frac{1}{\beta-\alpha}(\mathcal{M}_{1,0,(\alpha-1)}-\mathcal{M}_{1,0,(\beta-1)}),
\end{equation*}
\begin{equation*}\label{stmp}
  SP_{\alpha,\beta}=\frac{1}{\beta-\alpha}(\mathcal{M}_{1,(\alpha-1),0}-\mathcal{M}_{1,(\beta-1),0}),
\end{equation*}
\begin{equation*}\label{stmr}
  SRW_{\alpha,\beta}=\frac{1}{2(\beta-\alpha)}(\mathcal{M}_{2,0,(\alpha-1)}-\mathcal{M}_{2,0,(\beta-1)}),
\end{equation*}
\begin{equation*}\label{stmpw}
  SPW_{\alpha,\beta}=\frac{1}{2(\beta-\alpha)}(\mathcal{M}_{2,(\alpha-1),0}-\mathcal{M}_{2,(\beta-1),0}).
\end{equation*}


\noindent These representations in terms of PWM would enable us to develop inferential methods for all theses measures using the known results on PWM.
\vspace{-0.1in}
\section{Concluding remarks}\vspace{-0.1in}
The GMD is a well-known measure of dispersion that is used extensively in economics. We have established here some  relationships between the information measures and the GMD. We have shown that the GMD is a special case of the generalized cumulative residual entropy proposed recently by  Sudheesh et al. (2022). We have also established some relationships between the GMD and extropy measures. We have further presented some relationships between the dynamic versions of the cumulative extropy and the truncated GMD.

The probability weighted moments (PWM) generalize the concept of moments of a probability distribution. The estimates derived using probability weighted moments are often superior than the standard moment-based estimates. We have shown here that many of the information measures can be expressed in terms of PWM, which would enable us to develop inferential methods for these measures using the well-known results on PWM.

\end{document}